\documentclass{article}[12pt]
\usepackage{graphicx,pstricks,comment}
\usepackage[all]{xy}
\usepackage{subfigure}
\usepackage{rotating}
\usepackage{amsfonts}
\usepackage{amssymb}
\usepackage{amsmath}
\usepackage{amsthm}
﻿
\usepackage[utf8]{inputenc}
﻿
﻿
\usepackage{hyperref}
﻿
﻿
\usepackage{indentfirst}
﻿
\usepackage{epsfig}
\usepackage{psfrag}
\usepackage{enumerate}
\usepackage{array}

\makeatletter
\newtheorem*{rep@theorem}{\rep@title}
\newcommand{\newreptheorem}[2]{%
	\newenvironment{rep#1}[1]{%
		\def\rep@title{#2 \ref{##1}}%
		\begin{rep@theorem}}%
		{\end{rep@theorem}}}
\makeatother
﻿
﻿
﻿
\newcommand{\C}{\mathbb{C}}
\newcommand{\R}{\mathbb{R}}

\newcommand{\Z}{\mathbb{Z}}

\newcommand{\Q}{\mathbb{Q}}
\newcommand{\K}{\mathbb{K}}
\newcommand{\quat}{\mathbb{H}}

﻿
\newcommand{\HCd}{{\rm H}_\C^2}
\newcommand{\HCn}{{\rm H}_\C^n}

\newcommand{\la}{\langle}
\newcommand{\ra}{\rangle}

\newcommand{\onetwothreefour}{{Let $N$ be any Nil 3-manifold in families (1), (2),  (3) or (4) (Nil-tori, vertical or horizontal half-twist, or double half-twist). Then $N$ occurs as a cusp cross-section in the commensurability class of the Picard modular group ${\rm PU}(2,1,\mathcal{O}_d)$ for any squarefree $d \geqslant 1$.}}

\newcommand{\fivesixseven}{{The Nil 3-manifolds in families (5), (6), and (7) (1/4-twist, 1/3-twist and 1/6-twist) each occur as cusp cross-sections in a single arithmetic commensurability class.  For the 1/4-twist, this is the commensurability class of the Gauss-Picard modular group ${\rm PU}(2,1,\mathcal{O}_1)$, and for the 1/3-twist and 1/6-twist, this is the commensurability class of the Eisenstein-Picard modular group ${\rm PU}(2,1,\mathcal{O}_3)$.}}

\newcommand{\threefour}{{Let $N$ be any Nil 3-manifold in families (3) or (4) (horizontal half-twist or double half-twist). Then $N$ occurs as a cusp cross-section of an ${\rm H}^2_\C$-manifold, but not of any complex hyperbolic manifold or orbifold. That is, the image of any holonomy representation $\rho: \pi_1(N) \longrightarrow G_\infty \simeq {\rm Isom}({\rm Nil})$ must contain anitholomorphic isometries.}}

\newcommand{\nonarith}{{Every Nil 3-manifold occurs as the cusp cross-section of a non-arithmetic complex hyperbolic 2-manifold or ${\rm H}^2_\C$-manifold.}}
﻿

﻿
\newtheorem{conj}{Conjecture}
\newtheorem{thm}{Theorem}
\newreptheorem{thm}{Theorem}
\newtheorem{cor}{Corollary}%
\newtheorem{lem}{Lemma}
\newtheorem{prop}{Proposition}
\newtheorem{dfn}{Definition}
﻿\newtheorem{question}{Question}%

\newcommand{\Pf}{{\em Proof}. }
\newcommand{\EPf}{\hfill$\Box$\vspace{.5cm}}
﻿
﻿
﻿
\title{Nil 3-manifolds and cusps of complex hyperbolic surfaces}
\author{Julien Paupert, Connor Sell}
\date{\today}
\begin{document}
	
	\maketitle
	﻿
	\begin{abstract}
	McReynolds showed that every compact Nil 3-manifold occurs as the cusp cross-section of some arithmetic complex hyperbolic 2-manifold. We classify which commensurability classes of cusped, arithmetic, complex hyperbolic 2-manifolds admit cusps with cross-section homeomorphic to a given compact Nil 3-manifold.  In particular, there are some Nil 3-manifolds which occur as cusps in every such commensurability class, and some which only occur in a single commensurability class. We also show that every compact Nil 3-manifold occurs as the cusp cross-section of some non-arithmetic complex hyperbolic 2-manifold. 
	\end{abstract}

	\section{Introduction}
	
	A classical and well-studied question in topology is whether a given $ n $-manifold $ M $ occurs as the boundary of an $ (n+1) $-manifold $ W $.  For instance, Rohlin showed in \cite{R} that every compact orientable 3-manifold bounds a compact orientable 4-manifold - in other words, the cobordism group of 3-manifolds is trivial.  
	
	One can refine the question of bounding when the manifolds are equipped with geometric structures. For example, when $M$ is a compact hyperbolic $n$-manifold, one can ask whether $ M $ occurs as the totally geodesic boundary of a hyperbolic $(n+1)$-manifold-with-boundary $ W $. Analogously, when $M$ is a compact flat $n$-manifold one can ask whether $M$ occurs as the cusp cross-section of a (complete, finite-volume, 1-cusped) hyperbolic manifold $W$. (Recall that the ends of a complete, finite-volume, hyperbolic manifold, called \emph{cusps}, are diffeomorphic to $M \times \R^+$ with $M$ a compact flat $n$-manifold, called  a \emph{cusp cross-section}). In either case we will say that $M$ \emph{bounds geometrically}. Note that in the latter case, $M$ is diffeomorphic to the boundary of a manifold-with-boundary $W'$ such that $W=W' \setminus \partial W'$ admits a complete, finite-volume hyperbolic metric. ($W'$ is obtained from $W$ by adding the boundary $M$ ``at infinity").
	
	The latter version of geometric bounding also makes sense for the other negatively curved locally symmetric spaces, which do not admit any totally geodesic hypersurfaces (see Proposition 2.5.1 of \cite{CG}). These ``other" locally symmetric spaces - negatively curved but not constantly curved - are quotients of the hyperbolic spaces $ {\rm H}_\mathbb{K}^n $ with $\mathbb{K} = \mathbb{C} $ and $n \geqslant 2$,  $ \mathbb{K} =\mathbb{H} $ and $n \geqslant 2$, or $\mathbb{K} = \mathbb{O}$ and $n=2$. For these other locally symmetric spaces, the cusp cross-sections $M$ are now \emph{infra-nilmanifolds}, that is, compact quotients of a simply-connected nilpotent Lie group $\frak{N}$ by a discrete group of isometries. Gromov showed that the latter are exactly the compact \emph{almost flat manifolds} (see \cite{Gr}), so we will use the terms interchangeably. 
	
	More specifically, when $W$ is a cusped, negatively curved, locally symmetric space, its cusp cross-sections are quotients of 2-step nilpotent Lie groups usually called \emph{generalized Heisenberg groups} (over $\C$ or $\quat$, see e.g. \cite{CG}). For example, in the case of $\K=\C$ and $n=2$, which is the most relevant for this paper, $\frak{N}$ is the usual 3-dimensional Heisenberg group, which we will denote Nil, following the standard notation in 3-dimensional topology (this is one of the 8 Thurston geometries; see \cite{Sc} or Section~\ref{nil} for more details). An interesting feature of almost-flat manifolds is that, unlike the flat ones, they occur in infinite families in each dimension at least 3; the 3-dimensional ones are listed in Section~\ref{nil} by fundamental group.
		 	
	 
	 Farrell and Zdravkovska conjectured the following, Conjecture 2 of \cite{FZ} (in both statements $M$ is assumed to be compact and connected):
	 
	 \begin{conj}
	\begin{enumerate}
	 \item[(a)] If $M^n$ is a flat Riemannian manifold, then $M^n = \partial W^{n+1}$, where $W \setminus \partial W$ supports a complete hyperbolic structure with finite volume.
	 \item[(b)] If $M^n$ is an almost-flat Riemannian manifold, then $M^n = \partial W^{n+1}$, where $W \setminus \partial W$ supports a complete Riemannian metric with finite volume and negative sectional curvatures.
	 \end{enumerate}
	 \end{conj}
	 
	 Long and Reid showed in \cite{LR1} that (a) is false when $n=3$, providing examples of flat 3-manifolds which are not homeomorphic to any cross-section of a one-cusped (complete, finite-volume) hyperbolic 4-manifold. On the other hand, they also showed in \cite{LR2} that for any $n \geqslant 2$, every compact, flat $ n $-manifold occurs as the cusp cross-section of some hyperbolic $ (n+1) $-orbifold \cite{LR2} (possibly with several cusps).  This result was later upgraded by McReynolds in  \cite{M2} to say that such a flat manifold always occurs as the cusp cross-section of a hyperbolic $(n+1) $-\emph{manifold}, also possibly with several cusps.
	 
The restriction to manifolds (or even orbifolds) with a single cusp seems difficult, to say the least. In fact, while there is an abundance of 1-cusped hyperbolic surfaces and 3-manifolds (e.g. knot complements), it was only in 2013 that Kolpakov--Martelli produced in \cite{KM} the first example of a one-cusped hyperbolic 4-manifold. As far as we know, it is still an open question whether or not there exists a one-cusped hyperbolic $n$-manifold for $n \geqslant 5$. Analogously, it was not until very recently that Deraux (\cite{Der}) and Deraux--Stover (\cite{DS}) found the first examples of one-cusped complex hyperbolic 2-manifolds.

Moreover, it may well be that there do not exist any one-cusped hyperbolic or complex hyperbolic manifolds, or even orbifolds, in high enough dimensions. In the arithmetic world, Stover \cite{St1} showed that no arithmetic hyperbolic orbifold of dimension at least 30 can have exactly one cusp, and that similar bounds exist for arithmetic complex and quaternion hyperbolic orbifolds as well. (Since finite-volume quaternion hyperbolic orbifolds of dimension at least 2 are all arithmetic, the result is in fact unconditional in that case). Stover also produced in \cite{St1} examples of one-cusped (arithmetic) hyperbolic orbifolds in dimensions $n=10$ and 11; such orbifolds were known to exist for $n \leqslant 9$ among hyperbolic Coxeter groups. 

 
In view of these results it seems reasonable to relax the Farrell--Zdravkovska Conjecture to ask the following (where, again, the manifolds $M$ are assumed to be compact and connected):

\begin{question}
\begin{enumerate}
\item[(a)] Given a flat $n$-manifold $M$, is $M$ diffeomorphic to a cusp cross-section of some finite-volume hyperbolic $(n+1)$-manifold $W$? If so, of which such manifolds $W$?
\item[(b)] Given an almost-flat $n$-manifold $M$, is $M$ diffeomorphic to a cusp cross-section of some finite-volume complex, quaternion or octonion-hyperbolic manifold $W$ of real dimension $n+1$? If so, of which such manifolds $W$? 
\end{enumerate}
\end{question}	 
	 
Note that the class of complex, quaternion or octonion-hyperbolic $(n+1)$-manifold $W$ appearing in part (b) is, in case $W$ has a single cusp, exactly the class of locally symmetric Riemannian manifolds $W$ such that  $M = \partial W$, and $W \setminus \partial W$ supports a complete, finite-volume Riemannian metric with negative sectional curvatures, as in part (b) of the Farrell--Zdravkovska Conjecture. (Recall that a complex, quaternion or octonion-hyperbolic manifold is an isometric quotient of hyperbolic space $ {\rm H}_\K^n $ for $ \K = \mathbb{C}, \mathbb{H}$ or $\mathbb{O}$, see \cite{CG} and Section~\ref{complex} when $\K=\C$ for more details).	 

McReynolds studied in \cite{M1} the cusps of arithmetic real, complex, and quaternion-hyperbolic orbifolds, On one hand, he showed that every compact Nil 3-manifold occurs as the cusp cross-section of some (arithmetic) complex hyperbolic 2-orbifold (see below for some more precise terminology). This gives a positive answer to the first part of Question (b) above for any (compact, connected) almost-flat 3-manifold.

On the other hand, he showed that in all other cases (that is, complex dimension at least 3 and quaternion dimension at least 2) there exist (infinite families of) almost-flat manifolds modelled on the appropriate generalized Heisenberg group which are not diffeomorphic to a cusp cross-section of any finite-arithmetic complex or quaternion-hyperbolic manifold. As above, in the quaternion case all finite-volume manifolds are arithmetic and so this provides a negative answer to Question (b) above for infinitely many almost-flat manifolds in each dimension $4n+3$ with $n \geqslant 1$. 

Kamishima showed in \cite{K} that the cusp cross-section of any \emph{one-cusped} complex hyperbolic 2-manifold must have trivial holonomy (rotational part), that is, in the terminology of Section~\ref{nil} below, it must belong to the family (1) of Nil-tori. Moreover, by Corollary 2.2 of \cite{DS} the \emph{Euler number} (denoted $k$ in Section~\ref{nil}) of the Nil-torus must be a multiple of 4.
 The one-cusped complex hyperbolic 2-manifolds constructed in \cite{DS} have as cusp cross-sections Nil-tori with Euler number $k$ any odd multiple of 12; so in the above terminology such Nil-tori bound geometrically. 
	
	Both \cite{LR2} and \cite{M1} utilized arithmetic methods (and separability arguments) to construct manifolds with a given cross-section.  In both the real and complex cases, all known examples of hyperbolic manifolds in high enough dimensions are either arithmetic or use arithmetic manifolds in their constructions. (Here, high enough means dimension at least 4 in the complex case, and at least 6 in the real case).  
	
	The second author studied in \cite{Se} the second part of Question (a) for flat 3-manifolds. More precisely, the main results of \cite{Se} classify, given a flat 3-manifold $M$, in which commensurability classes of cusped, arithmetic hyperbolic 4-manifolds $M$ occurs as a cusp cross-section. In this paper, we study the analogous second part of Question (b) for almost-flat 3-manifolds, that is we study, given a Nil 3-manifold $M$, in which cusped, complex hyperbolic 2-manifolds $M$ can occur as a cusp cross-section. 
	
	In the arithmetic case, similarly to \cite{Se}, given a Nil 3-manifold $M$, we classify in which commensurability classes of cusped, arithmetic complex hyperbolic 2-manifolds $M$ occurs as a cusp cross-section. Such commensurability classes are enumerated by the \emph{Picard modular groups} ${\rm PU}(2,1,\mathcal{O}_d)$ with $d$ a positive square-free integer and $\mathcal{O}_d$ the ring of integers of $\Q(i\sqrt{d})$ (see Section~\ref{picard}).
	First, we show that some Nil 3-manifolds occur in every commensurability class of cusped arithmetic complex hyperbolic 2-manifolds (see Section~\ref{nil} for the labelling of families (1)--(7)).
	 
	\begin{thm}
		\label{1234}
		\onetwothreefour
	\end{thm}
	
	Secondly, we show that some Nil 3-manifolds are obstructed from some commensurability classes.  In fact, they each occur in only a single commensurability class of arithmetic complex hyperbolic 2-manifolds.  
	Note that every Nil 3-manifold is covered by either Theorem \ref{1234} or Theorem \ref{567}.
	
	\begin{thm}
		\label{567}
		\fivesixseven
	\end{thm}
	
	In passing, we prove the following related result.  This shows that certain Nil 3-manifolds $ N $ can only occur as cusp cross-sections in $ {\rm H}_\C^2 $-manifolds which are not complex hyperbolic manifolds (see Definition~\ref{manifolds}); in other words, any embedding of $ \pi_1(N) $ into $ \text{Isom}({\rm H}_\C^2) $ must contain antiholomorphic isometries in its image (see Section~\ref{complex}).
	
	\begin{thm}
		\label{34antihol}
		\threefour
	\end{thm}
	
	Finally, we combine results of \cite{M1}  and \cite{M2} with the detailed analysis of cusp groups in the known non-arithmetic complex hyperbolic orbifolds in \cite{DPP2} to show the following:
	
	\begin{thm}
		\label{NA}
		\nonarith
	\end{thm}
	
	The paper is organized as follows. Section \ref{sec:background} gives background on Nil-manifolds, complex hyperbolic geometry, Picard modular groups and their cusps and non-arithmetic complex hyperbolic 2-orbifolds and their cusps. Section \ref{sec:results} is devoted to the results stated above and their proofs.

	
	
	\section{Geometric background}
	\label{sec:background}
	﻿
	\subsection{Nil-manifolds}\label{nil}
	﻿
	The Heisenberg group Nil is a 3-dimensional Lie group which can be written as $\C \times \R$ with group law given by:
	\begin{equation}
		(z_1,t_1)\cdot (z_2,t_2)=(z_1+z_2,t_1+t_2+2{\rm Im}(z_1\overline{z_2})).
	\end{equation}
	It is a 2-step nilpotent group, as one can easily check that $[{\rm Nil},{\rm Nil}]=Z({\rm Nil})=\{ (0,t) \, | \, t \in \R \}$, the \emph{vertical axis}.
	This model of Nil naturally arises from the action of ${\rm PU(2,1)}$ on the 3-sphere $\partial_\infty {\rm H}^2_\C$ as we will see in the next subsection. (The most standard model of Nil is the subgroup of upper-triangular unipotent matrices in ${\rm SL}_3(\R)$, see \cite{Sc}). 
	
	﻿
	As any Lie group, Nil admits a left-invariant Riemannian metric, unique up to scale. By construction, (left-) translations are isometries for this metric. We will describe the isometry group ${\rm Isom}({\rm Nil})$ in more detail in the next subsection, but for now we just note that, as in the Euclidean case, isometries are \emph{affine diffeomorphisms} of Nil in the sense that:  
	$${\rm Isom}({\rm Nil})={\rm Nil} \rtimes C$$
	for some compact Lie subgroup $C$ of ${\rm Aut}({\rm Nil})$. In fact we will see below that $C={\rm U}(1) \rtimes \Z_2$, with the ${\rm U}(1)$ factor acting on Nil by rotations about the (vertical) $t$-axis, and $\Z_2$ acting by $(z,t) \mapsto (\overline{z},-t)$.
	﻿
	\begin{dfn} A \emph{Nil-manifold} is a compact Riemannian 3-manifold $N$ locally modelled on ${\rm Nil}$ (in the sense of Ehresmann--Thurston geometric structures). In other words, $N \simeq {\rm Nil}/\Delta$, where $\Delta \simeq \pi_1(N)$ is a discrete torsion-free subgroup of ${\rm Isom}({\rm Nil})$. 
	\end{dfn}
	﻿
	Note that there are various definitions of Nil- (or nil-)manifolds in the literature, but this is the standard definition from 3-manifold topology/geometric structures. (We will quote the classification from \cite{Dek} below, where they are called infra-nilmanifolds).
	﻿
	
	Given a Nil-manifold $N$, we will be interested in recovering the geometric structure on $N$ from representations of $\pi_1(N)$. To this end, recall that a representation $\rho : \pi_1(N) \longrightarrow {\rm Isom}({\rm Nil})$ is a \emph{holonomy representation} (for the Nil-structure on $N$) if ${\rm Nil}/\rho(\pi_1(N))$ is diffeomorphic to $N$. The basic theory of geometric structures (see e.g. \cite{Go2}, Corollary 8.1.4) gives us the following:
	
	\begin{lem}\label{holonomy} Let $N$ be a Nil-manifold and $\rho : \pi_1(N) \longrightarrow {\rm Isom}({\rm Nil})$ a discrete and faithful representation with ${\rm Nil}/\rho(\pi_1(N))$ compact. Then ${\rm Nil}/\rho(\pi_1(N))$ is diffeomorphic to $N$; that is, $\rho$ is a holonomy representation of $\pi_1(N)$. 
	\end{lem}
	
	As a special case of the second generalized Bieberbach theorem below, the holonomy representation of a Nil-manifold group $\rho : \pi_1(N) \longrightarrow {\rm Isom}({\rm Nil})$ is unique up to conjugation by the affine group ${\rm Nil} \rtimes {\rm Aut}({\rm Nil})$. The following generalized Bieberbach theorems are due to Auslander (\cite{A}), with a corrected statement and proof for the second statement by Lee--Raymond (\cite{LR}, see also discussion in Section 2.2 of \cite{Dek}).
	
	\begin{thm}\label{GenBieberbach} Let $N$ be a simply-connected nilpotent Lie group and $C$ a compact Lie subgroup of ${\rm Aut}(N)$. 
		\begin{enumerate}
			\item(\cite{A}) Let $\Gamma$ be a discrete cocompact subgroup of $G=N \rtimes C$. Then $\Gamma \cap N$ is a cocompact lattice in $N$ and $\Gamma/(\Gamma \cap N)$ is finite.
			\item (\cite{LR}) Let $\Gamma_1$,  $\Gamma_2$ be discrete cocompact subgroups of $G=N \rtimes C$. If $\Gamma_1$ and $\Gamma_2$ are isomorphic (as abstract groups) then they are conjugate in the affine group $N \rtimes {\rm Aut}(N)$.
		\end{enumerate}
	\end{thm}
	Note that the (first two) classical Bieberbach theorems are obtained by taking $N=\R^n$ (and $C={\rm O}(n)$) in these statements. The third Bieberbach theorem, stating that up to conjugacy, there are only finitely many such $\Gamma$ in each dimension, is no longer true for more general nilpotent groups. There is a meaningful generalization of this third Bieberbach theorem (see Section 2.3 of \cite{Dek}), but when $N={\rm Nil}$ there is in fact an explicit classification of Nil-manifold groups, coming in seven infinite families as follows. 
	﻿
	
	We give presentations for the fundamental groups of the manifolds in these seven families below; the numbering is the same as in \cite{M1} (these are the \emph{AB-groups} from \cite{Dek}, appearing on pp. 159--166 in families numbered 1, 2, 4, 8, 10, 13, 16 respectively). We give these families more evocative names for ease of reference. 
	
	\begin{enumerate}
		\item (Nil-tori) $ \langle a, b, c ~|~ [b,a] = c^k, [c,a] = [c,b] = 1 \rangle $ \\ with $ k \in \mathbb{N} $
		\item (Vertical half-twist) $ \langle a, b, c, \alpha ~|~ [b,a] = c^k, [c,a] = [c,b] = [c, \alpha] = 1, \alpha a = a^{-1} \alpha, \alpha b = b^{-1} \alpha, \alpha^2 = c \rangle $ \\ with $ k \in 2\mathbb{N} $.
		\item (Horizontal half-twist) $ \langle a, b, c, \alpha ~|~ [b,a] = c^{2k}, [c,a] = [c,b] = [a, \alpha] = 1, \alpha b = b^{-1} \alpha c^{-k}, \alpha c = c^{-1} \alpha, \alpha^2 = a \rangle $ \\ with $ k \in \mathbb{N} $.
		\item (Double half-twist) $ \langle a, b, c, \alpha, \beta ~|~ [b,a] = c^{2k}, [c,a] = [c,b] = [c, \alpha] = [a, \beta] = 1, \alpha a = a^{-1} \alpha c^k, \\ \alpha b = b^{-1} \alpha c^{-k}, \beta b = b^{-1} \beta c^{-k}, \beta c = c^{-1} \beta, \alpha \beta = a^{-1} b^{-1} \beta \alpha c^{-k-1}, \alpha^2 = c, \beta^2 = a \rangle $ \\ with $ k \in 2 \mathbb{N} $.
		\item ($1/4$-twist) $ \langle a, b, c, \alpha ~|~ [b,a] = c^k, [c,a] = [c,b] = [c, \alpha] = 1, \alpha a = b \alpha, \alpha b = a^{-1} \alpha, \alpha^4 = c^p \rangle $ \\ with $ k \in 2 \mathbb{N} $ and either $ p = 1 $ or $ k \in 4\mathbb{N} $ and $ p = 3 $.
		\item  ($1/3$-twist) $ \langle a, b, c, \alpha ~|~ [b,a] = c^k, [c,a] = [c,b] = [c, \alpha] = 1, \alpha a = b \alpha c^{k_1}, \alpha b = a^{-1} b^{-1} \alpha, \alpha^3 = c^{k_2} \rangle $ \\ with $ k > 0 $ and either $ k \equiv 0 $ (mod 3), $ k_1 = 0 $, and $ k_2 = 1 $ or 2, \\ or $ k \equiv 1,2 $ (mod 3), $ k_1 = 1 $, and $ k_2 = 1 $.
		\item  ($1/6$-twist) $ \langle a, b, c, \alpha ~|~ [b,a] = c^k, [c,a] = [c,b] = [c, \alpha] = 1, \alpha a = ab \alpha, \alpha b = a^{-1} \alpha, \alpha^6 = c^{k_1} \rangle $ \\ with $ k > 0 $ and either $ k \equiv 0,4 $ (mod 6) and $ k_1 = 1 $, or $ k \equiv 0, 2 $ (mod 6) and $ k_1 = 5 $.
	\end{enumerate}
	
	In the geometric realization of these groups (in the holonomy representation), the generators $ a $, $ b $, and $ c $ are Heisenberg translations,  with $a, b$ horizontal translations and $c$ a vertical translation (see next subsection for terminology).  The generators $ \alpha $ and $ \beta $ are ellipto-parabolic, and play a similar role to the generators of flat 3-manifolds with nontrivial rotational part.
	
	In this way, there is a correspondence between these seven families of Nil-manifolds and the six orientable flat 3-manifolds, with families (1), (4), (5), (6), and (7) corresponding to the 3-torus, Hantzsche-Wendt manifold, $1/4$-twist, $1/3$-twist, and $1/6$-twist, respectively.  Families (2) and (3) both correspond to the half-twist, with the twist occurring along a vertical or horizontal Heisenberg translation, respectively.  In fact, setting $ k = 0 $ in these presentations (along with particular values for $ p $, $ k_1 $, or $ k_2 $) yields the fundamental group of the corresponding flat manifold.
	
	\subsection{Complex hyperbolic geometry}\label{complex}
	
	See \cite{CG} or \cite{Go1} for a more
 detailed discussion of the background material in this section.
 \medskip
 
 	{\bf Projective models:} Consider $\C^{n,1}$, the vector space $\C^{n+1}$ endowed with a Hermitian form $\langle \cdot \, , \cdot \rangle$ of signature $(n,1)$.
	Let $V^-=\left\lbrace Z \in \C^{n,1} | \langle Z , Z \rangle <0 \right\rbrace$ and $\pi: \C^{n+1}-\{0\} \longrightarrow \C{\rm P}^n$ denote projectivization.
	Define ${\rm H}_\C^n$ to be $\pi(V^-) \subset \C{\rm P}^n$, endowed with the distance $d$ given by:
	﻿
	\begin{equation}\label{dist}
		\cosh ^2 \frac{1}{2}d(\pi(X),\pi(Y)) = \frac{|\langle X, Y \rangle|^2}{\langle X, X \rangle  \langle Y, Y \rangle}
	\end{equation}
	﻿
	This distance is in fact induced by a complete Riemannian metric (called the \emph{Bergman metric}) which has negative sectional curvatures $1/4$-pinched (between $-1$ and $-1/4$ in the above normalization).
	
	Different choices of Hermitian forms of signature $(n,1)$ give rise to different models of $\HCn$. For our purposes we will consider the  \emph{Siegel model}, generalizing the Poincar\'e upper half-plane model of ${\rm H}^2_\R$, where the Hermitian form is given by $\la Z,Z\ra=2{\rm Re} (z_1\overline{z_{n+1}}) +|z_2|^2+\cdots +|z_n|^2$. This form provides especially nice coordinates for the stabilizer of a particular point at infinity, in the form of upper-triangular matrices. 
	﻿
	\medskip
	
	{\bf Isometries:} From \eqref{dist} it is clear that ${\rm PU}(n,1)$ acts by isometries
	on ${\rm H}_\C^n$, where ${\rm U}(n,1)$ denotes the subgroup of ${\rm
		GL}(n+1,\C)$ preserving $\langle \cdot , \cdot \rangle$, and ${\rm
		PU}(n,1)$ its image in ${\rm PGL}(n+1,\C)$. In fact, PU($n$,1) is the group of holomorphic isometries of 
	${\rm H}_\C^n$, and the full group of isometries ${\rm Isom}({\rm H}_\C^n)$ is ${\rm PU}(n,1) \rtimes \Z/2$, where the $\Z/2$ factor 
	corresponds to a \emph{real reflection}, that is, an antiholomorphic involution (see below). 
	﻿
	\begin{dfn}\label{manifolds} A \emph{complex hyperbolic manifold} is a quotient ${\rm H}_\C^n/\Gamma$ where $\Gamma$ is a torsion-free discrete subgroup of ${\rm PU}(n,1)$. More generally, a \emph{manifold modelled on ${\rm H}_\C^n$}, or in short, \emph{${\rm H}_\C^n$-manifold}, is a quotient  ${\rm H}_\C^n/\Gamma$ with $\Gamma$ a torsion-free discrete subgroup of ${\rm Isom}({\rm H}_\C^n)$.
	\end{dfn}
	
	Note that, since  ${\rm PU}(n,1)$ consists of holomorphic isometries, the complex (in fact K\"ahler) structure of ${\rm H}_\C^n$ descends to complex hyperbolic manifolds. On the other hand, if $\Gamma$ is a torsion-free discrete subgroup of ${\rm Isom}({\rm H}_\C^n)$ containing an antiholomorphic isometry then the complex structure does not descend to the quotient ${\rm H}_\C^n/\Gamma$, and this quotient in fact does not admit any complex structure. 
	﻿

	﻿
	Holomorphic isometries of $\HCn$ can be of three types, depending on the number 
	and location of their fixed points. Namely, $g \in {\rm PU}(n,1) - \{ {\rm Id} \}$ is :
	\begin{itemize}
		\item \emph{elliptic} if it has a fixed point in ${\rm H}_\C^n$
		\item \emph{parabolic} if it has (no fixed point in ${\rm H}_\C^n$ and)
		exactly one fixed point in $\partial{\rm H}_\C^n$
		\item \emph{loxodromic}: if it has (no fixed point in ${\rm H}_\C^n$ and) exactly two fixed points in $\partial{\rm H}_\C^n$
	\end{itemize}
	﻿
	A parabolic isometry is called {\it unipotent} (or a \emph{Heisenberg translation}) if it has a unipotent lift in ${\rm
		U}(n,1)$; if not it is called {\it ellipto-parabolic}.  In dimensions $n>1$, unipotent 
	isometries are either  {\it 2-step} or {\it 3-step}, according to whether the minimal polynomial of their unipotent lift
	is $(X-1)^2$ or $(X-1)^3$ (see section 3.4 of \cite{CG}). We will usually call 2-step nilpotent isometries \emph{vertical Heisenberg translations} and 3-step nilpotent isometries \emph{horizontal Heisenberg translations}. This terminology reflects the position of the subspaces preserved by each type of isometry	in the Heisenberg group, relative to its standard contact structure. Any two horizontal (resp. vertical) Heisenberg translations are conjugate under an isometry.
	
	{\bf Totally geodesic subspaces and related isometries:} 
	A {\it complex k-plane} is a projective $k$-dimensional subspace of 
	$\C P^n$ intersecting $\pi(V^-)$
	non-trivially (so, it is an isometrically embedded copy of ${\rm
		H}_\C^{k} \subset {\rm H}_\C^n$). Complex 1-planes are usually
	called {\it complex lines}. 
	
	A {\it real k-plane} is the projective image of a totally real
	$(k+1)$-subspace $W$ of $\C^{n,1}$, i. e. a $(k+1)$-dimensional real
	vector subspace such that $\langle v,w \rangle \in \R$ for all $v,w
	\in W$.  We will usually call real 2-planes simply real planes, or
	$\R$-planes. Every real $n$-plane in ${\rm H}_\C^n$ is the
	fixed-point set of an antiholomorphic isometry of order 2 called a {\it real
		reflection} or $\R$-reflection. The prototype of such an isometry
	is the map given in affine coordinates by:
	\begin{equation}\label{sigma}
		\sigma: (z_1,...,z_n) \mapsto
		(\overline{z_1},...,\overline{z_n})
	\end{equation}
	This is an isometry provided that the Hermitian form has real coefficients.
	﻿
	﻿
	\medskip
	
	{\bf The Siegel model in dimension 2:}
	As mentioned in the previous section, this 
	model corresponds to the Hermitian form given by the matrix:
	\begin{equation}\label{Hermitian}
		H=
		\begin{bmatrix}
			0 & 0 & 1 \\
			0 & 1 & 0 \\
			1 & 0 & 0 
		\end{bmatrix}
	\end{equation}
	In this model, any point  $p \in \HCd$ admits a unique lift to $\C^{2,1}$ of the following form, called its \emph{standard lift}:
	\begin{equation}
		{\bf p}=\begin{bmatrix}
			(-|z|^2-u+it)/2\\z\\1
		\end{bmatrix}\mbox{ with } (z,t,u)\in\C\times\R\times (0,\infty) .
	\end{equation}
	﻿
	\noindent The coordinates $(z,t,u)$ are called \textit{horospherical coordinates} of $p$. The boundary at infinity $\partial_\infty\HCd$ is 
	the level set $\lbrace u=0\rbrace$, together with the distinguished point at infinity, given by
	﻿
	$$
	q_\infty\sim\left[
	\begin{array}{c}
		1 \\ 0 \\ 0
	\end{array}\right].
	$$
	
	Level sets $\lbrace u=u_0\rbrace$ with fixed $u_0>0$ are called \emph{horospheres based at $q_\infty$}. 
	The punctured boundary at infinity $\partial_\infty\HCd\setminus \{ q_\infty \}$ is a copy of the Heisenberg group ${\rm Nil}$, identified as above with $\C \times \R$ with group law given by:
	﻿
	\begin{equation}
		\label{Heisprod}(z_1,t_1)\cdot (z_2,t_2)=(z_1+z_2,t_1+t_2+2{\rm Im}(z_1\overline{z_2})).
	\end{equation}
	﻿
	The stabilizer of $q_\infty$ in ${\rm PU}(2,1)$ corresponds to upper triangular matrices, and is generated by the following 3 types of isometries: Heisenberg translations $T_{(z,t)}$ ($(z,t)\in \C \times \R$), Heisenberg rotations $R_u$ ($u \in {\rm U}(1)$) and Heisenberg dilations $D_r$ ($r >0$), where:
	﻿
	\begin{equation}\label{stabinf}
		\begin{array}{ccc}
			T_{(z,t)}=
			\begin{bmatrix}
				1 & -\overline{z} & -(|z|^2-it)/2 \\
				0 & 1 & z \\
				0 & 0 & 1
			\end{bmatrix}
			&
			R_u=
			\begin{bmatrix}
				1 & 0 & 0 \\
				0 & u & 0 \\
				0 & 0 & 1
			\end{bmatrix}
			&
			D_r=
			\begin{bmatrix}
				r & 0 & 0 \\
				0 & 1 & 0 \\
				0 & 0 & 1/r
			\end{bmatrix}.
		\end{array}
	\end{equation}
	﻿
	In Heisenberg coordinates, these correspond to the following:
	\begin{itemize}
		\item $T_{(z,t)}$ is left multiplication by $(z,t)$: $(w,s)\longmapsto (z,t)\cdot (w,s)$,
		\item $R_u$ is given by $(w,s)\longmapsto (uw,s)$,
		\item $D_r$ is the Heisenberg dilation $(w,s)\longmapsto (rw,r^2s)$.
	\end{itemize}
	﻿
	Heisenberg translations and rotations preserve each horosphere based at $q_\infty$ whereas
	Heisenberg dilations permute horospheres based at $q_\infty$. 
	It turns out that the subgroup of ${\rm Stab}_{{\rm PU}(2,1)}(q_\infty)$ preserving each horosphere based at $q_\infty$ is ${\rm Isom}^0({\rm Nil})$, the identity component of isometries for the Riemannian metric on Nil discussed in the previous section. 
	﻿
	Note that ${\rm Isom}^0({\rm Nil})= \{ P_{(z,t,u)} \, | \, z \in \C, t \in \R, u \in {\rm U}(1) \}$, where:
	﻿
	\begin{equation}\label{Pztu}
		P_{(z,t,u)}=T_{[z,t]}R_u =\begin{bmatrix}
			1 & -u\overline{z} & -(|z|^2-it)/2\\
			0 & u & z\\
			0 & 0 & 1       
		\end{bmatrix} 
	\end{equation}
	Moreover, denoting as above $\sigma$ complex conjugation in affine coordinates of ${\rm H}^2_\C$, $\sigma$ acts on Nil by $(z,t) \mapsto (\overline{z},-t)$ and:
	﻿
	$${\rm Isom}({\rm Nil})={\rm Isom}^0({\rm Nil}) \cup {\rm Isom}^0({\rm Nil}) \sigma$$ .
	
	For notational convenience we will denote $G={\rm Isom}({\rm H}^2_\C)=\langle {\rm PU}(2,1), \sigma \rangle$, $G^0={\rm Isom}^0({\rm H}^2_\C)={\rm PU}(2,1)$ the subgroup of holomorphic isometries, and $G_\infty$, $G^0_\infty$ their respective subgroups preserving ($q_\infty$ and) each horosphere based at $q_\infty$. That is:
	\begin{eqnarray}\label{ParabolicSubgroups}
		& & G^0_\infty={\rm Isom}^0({\rm Nil})= \left\lbrace P_{(z,t,u)} \, | \, z \in \C, t \in \R, u \in {\rm U}(1) \right\rbrace \\
		& & G_\infty={\rm Isom}({\rm Nil})= \langle G^0_\infty, \sigma \rangle = G^0_\infty \cup G^0_\infty \sigma. 
	\end{eqnarray}
	﻿
	 For any subring $R$ of $\C$ we likewise denote $G^0(R)$ (resp. $G_\infty^0(R)$) the subgroup of $G^0$ (resp. $G^0_\infty$) consisting of (images in PU(2,1) of) matrices with entries in $R$. In the next subsection we will consider the groups with $R={\mathcal O}_d$, the ring of integers of $\Q(i\sqrt{d})$, giving the Picard modular groups $G^0({\mathcal O}_d)={\rm PU}(2,1,{\mathcal O}_d)$ and the corresponding cusp subgroup $G^0_\infty{(\mathcal O}_d)$. By extension we will also denote $G({\mathcal O}_d)=\langle G^0({\mathcal O}_d), \sigma \rangle$ and $G_\infty({\mathcal O}_d)=\langle G^0_\infty({\mathcal O}_d), \sigma \rangle$ their index-2 extensions obtained by adding the $\R$-reflection $\sigma$ defined in (\ref{sigma}).
	   
	   \medskip
	   
	A (holomorphic) parabolic isometry is conjugate in ${\rm Isom}(\HCd)$ to exactly one of the following:
	\begin{itemize}
		\item $P_{(1,0,0)}=T_{[1,0]}$ if it is a horizontal Heisenberg translation (3-step unipotent),
		\item $P_{(0,1,0)}$ if it is a vertical Heisenberg translation (2-step unipotent),
		\item $P_{(0,1,\theta)}$ for some non-zero $\theta\in\R/2\pi\Z$ if it is ellipto-parabolic. 
	\end{itemize}
	﻿
	﻿
	Also note that Heisenberg dilations act by conjugation as: $D_r P_{(z,t,u)} D_r^{-1}=P_{(rz,r^2t,u)}$ (for any $z \in \Z, t\in \R, u \in {\rm U}(1), r>0$).
	﻿
	Recall also from \cite{Sc} that the exact sequence:
	\begin{equation}\label{Pi} 1 \longrightarrow \R \longrightarrow {\rm Nil} \overset{\Pi}{\longrightarrow} \C \longrightarrow 1
	\end{equation}
	induces an exact sequence:
	\begin{equation}\label{exact}
		1 \longrightarrow{\rm Isom} (\R) \longrightarrow {\rm Isom} ({\rm Nil}) \overset{\Pi_*}{\longrightarrow} {\rm Isom} (\C) \longrightarrow 1.
	\end{equation}
	Explicitly: $\Pi_*(P_{(z,t,\theta)})=\left(\begin{matrix} e^{i\theta} & z \\ 0 & 1\end{matrix}\right)$, acting on $\C$ by 
	$w \mapsto e^{i\theta}w+z$, and $\Pi_*(\sigma)$ is complex conjugation in $\C$.
	﻿
	﻿
	﻿
	In practice, we will write elements of $G^0_\infty$ as triples $(z,t,u)$ and compute products in $G_\infty$ by the following rules, which are straightforward from the matrix form $P_{(z,t,u)}$ above.
	﻿
	\begin{lem}\label{ParabolicProducts} Let $z_1,z_2,z \in \C$, $t_1,t_2,t \in \R$ and $u_1,u_2,u \in {\rm U}(1)$. Then:
		\begin{enumerate}
			\item $(z_1,t_1,u_1)\cdot (z_2,t_2,u_2) = (u_1z_2+z_1,t_1+t_2+2{\rm Im}(\overline{u_1}z_1\overline{z_2}), u_1u_2)$,
			\item $\sigma (z,t,u) = (\overline{z}, -t, \overline{u}) \sigma$.
		\end{enumerate}
	\end{lem} 
	In some of the representations below we will consider elements whose square is a vertical Heisenberg translation $(0,t,1)$ (for some $t \in \R$). For this it is instructive to note that for any $z \in \C, t \in \R$, $(z,t,-1)^2=(0,2t,1)$, that is: the square of any ellipto-parabolic with rotational part $-1$ is a vertical Heisenberg translation. At the other extreme:
	﻿
	\begin{lem} No antiholomorphic isometry in $G_\infty$ squares to a vertical Heisenberg translation. That is, for any $(z,t,u) \in G^0_\infty$, if $((z,t,u)\sigma)^2=(0,t',1)$ for some $t' \in \R$ then $t'=0$. 
	\end{lem}
	﻿
	\Pf By the above lemma, $((z,t,u)\sigma)^2=(z,t,u) \cdot(\overline{z}, -t, \overline{u})=(u\overline{z}+z	, -2{\rm Im}(u\overline{z}^2),1)$. Now if $u\overline{z}+z=0$, $u\overline{z}^2=-z\overline{z}$ and so ${\rm Im}(u\overline{z}^2)=0$. \EPf
	﻿
	﻿
	\subsection{Picard modular groups and their cusps}\label{picard}
	
	Although the Nil-manifolds are more complicated than the flat 3-manifolds, the commensurability classes of cusped arithmetic complex hyperbolic 2-manifolds are easier to understand than that of their real counterparts.
	
	A group of the form $ {\rm PU}(2, 1, \mathcal{O}_d) $, where $ d \geqslant 1$ is a squarefree integer and $ \mathcal{O}_d $ is the ring of integers of $ \mathbb{Q}(i\sqrt{d}) $, is called a \emph{Picard modular group}.  (Here we think of  ${\rm PU}(2, 1)$ as corresponding to the Siegel model, with Hermitian form  given in (\ref{Hermitian}). See Section 3.1 of \cite{St2} for a more general discussion involving different forms.) 
	
	Every non-cocompact, arithmetic lattice in $ {\rm PU}(2, 1, \mathcal{O}_d) $ is commensurable to a Picard modular group ${\rm PU}(2,1, \mathcal{O}_d) $ for some $ d $, and no two distinct ${\rm PU}(2, 1, \mathcal{O}_d)$ are commensurable (see \cite{St2}). In other words, passing to the quotient orbifolds $ {\rm H}_\C^2 / {\rm PU}(2,1, \mathcal{O}_d) $, the commensurability classes of arithmetic complex hyperbolic 2-manifolds are in bijective correspondence with the Picard modular groups. (Note that, as all Picard modular groups contain torsion (see below), the corresponding quotients are not manifolds.)
	﻿
	\medskip
	
	{\bf Cusp subgroups of Picard modular groups:}
	
	\medskip
	﻿
	The Picard modular groups are non-cocompact lattices, hence contain parabolic elements. In fact, it is clear from the matrix form (\ref{Pztu}) that they all have the point $q_\infty=[1,0,0]^T$ as a cusp point (that is, a fixed point of a parabolic isometry in the group).
	We will denote by $\Gamma_\infty(d) $ the corresponding cusp subgroup of $ {\rm PU}(2, 1, \mathcal{O}_d) $, that is, the stabilizer of $q_\infty$ in $ {\rm PU}(2, 1, \mathcal{O}_d) $. (In the above notation \ref{ParabolicSubgroups}, $\Gamma_\infty(d) =G^0\infty({\mathcal O}_d)$). 
	
	These subgroups have been explicitly described in \cite{FP} ($d=3$), \cite{FFP} ($d=1$) and \cite{PW} (all other $d$), in terms of the exact sequence (\ref{exact}). We will denote by  $\Gamma_\infty^v(d) $, $\Gamma_\infty^h(d)$ the \emph{vertical} and \emph{horizontal} parts of $\Gamma_\infty(d) $, that is the kernel and image of the homomorphism $\Pi_*$, so that the exact sequence becomes, when restricted to $\Gamma_\infty(d) $:
	﻿
	\begin{equation}\label{exactd}
		1 \longrightarrow \Gamma_\infty^v(d) \longrightarrow \Gamma_\infty(d) \overset{\Pi_*}{\longrightarrow} \Gamma_\infty^h(d)  \longrightarrow 1.
	\end{equation}
	﻿
	Combining Proposition~3.1 of \cite{FP}, Proposition~2 of \cite{FFP} and Lemma~5.2 of \cite{PW} gives the following, where we denote $\Delta(p,q,r)$ the orientation-preserving $(p,q,r)$-triangle group in ${\rm Isom}(\C)$:
	﻿
	\begin{prop}\label{PicardCusp} 
		\begin{itemize}
			\item[(a)] For all $d$, $\Gamma_\infty^v(d)$ is generated by the vertical Heisenberg translation $T_{(0,2\sqrt{d})}$. 
			\item[(b)] $\Gamma_\infty^h(3)=\Delta(2,3,6)$; $\Gamma_\infty^h(1)$ has index 2 in ${\rm Isom}^+(\Z[i])=\Delta(2,4,4)$. 
			\item[(c)]  For all other $d$, $\Gamma_\infty^h(d) = {\rm Isom}^+(\mathcal{O}_d)$ if $d \equiv 3 \ (mod \, 4)$;   $\Gamma_\infty^h(d)$ has index 2 in ${\rm Isom}^+(\mathcal{O}_d)$ if $d \equiv 1,2 \ (mod \, 4)$.
		\end{itemize}
	\end{prop}
	
	{\bf Remarks:} 
	\begin{itemize}
		\item In part (a), Lemma~5.2 of \cite{PW} mistakenly stated $T_{(0,\sqrt{d})}$ instead of $T_{(0,2\sqrt{d})}$ (the former has entry $i\sqrt{d}/2$ in the top-right corner, which is not in $\mathcal{O}_d$).
		\item Recall that $\mathcal{O}_d=\Z[\tau]$, where $\tau=i\sqrt{d}$ if $d \equiv 1,2 \ (mod \, 4)$ and $\tau=\frac{1+i\sqrt{d}}{2}$ if $d \equiv 3 \ (mod \, 4)$. In part (c), for $d\neq 1,3$, ${\rm Isom}^+(\mathcal{O}_d)$ is generated by the translations by $1$ and $\tau$ and a half-turn. The corresponding half-turn is a torsion element in $\Gamma_\infty(d)$ for all $d$.
	\end{itemize}
	﻿
	\subsection{Non-arithmetic lattices in ${\rm PU}(2,1)$ and their cusps}
	﻿
	Non-arithmetic lattices in PU(2,1) were constructed by Mostow and Deligne--Mostow in the 1980's, and more recently by Deraux, Parker and the first author in \cite{DPP1}, \cite{DPP2}. The latter paper classified the known examples up to commensurability, showing that they fall into 22 commensurability classes, 10 of which consist of cusped lattices (Table~11 of \cite{DPP2} gives 11 such classes, but the class containing $\Gamma(3,1/12)$ and $\Gamma(8,7/24)$ is in fact cocompact, see Appendix A.9). A key geometric fact is that each of the 22 classes contains a group generated by 3 complex reflections with the same order $p \geqslant 3$. 
	﻿
	
	We do not need to review details of the construction or notation; we simply collect here the information about the cusp groups for those lattices that are non-cocompact, see Section~7.2.2 of \cite{DPP2}. These all have $p=3,4$ or 6, and the cusp group $\Gamma_\infty$ is generated by two complex reflections $A,B$ of order $p$. As we did for the Picard modular groups in (\ref{exactd}), we characterize each of these cusp groups by its horizontal part  $\Gamma_\infty^h$ and vertical part $\Gamma_\infty^v$ in the exact sequence (\ref{exact}). In all three cases, $\Gamma_\infty^h$ is a triangle group and $\Gamma_\infty^v$ is generated by a vertical translation $T$ which is a power of $AB$. More precisely:
	﻿
	\begin{lem}
		\begin{itemize}
			\item When $p=3$, $\Gamma_\infty^h=\Delta(3,3,3)$ and $\Gamma_\infty^v=\la T \ra$ with $T=(AB)^3$.
			\item When $p=4$, $\Gamma_\infty^h=\Delta(2,4,4)$ and $\Gamma_\infty^v=\la T \ra$ with $T=(AB)^2$.
			\item When $p=6$, $\Gamma_\infty^h=\Delta(2,3,6)$ and $\Gamma_\infty^v=\la T \ra$ with $T=(AB)^3$.
		\end{itemize}
	\end{lem}
	For completeness we note that all three cases above occur as a cusp group for some non-arithmetic lattice, see again Table~11 of \cite{DPP2}.
	Comparing this with the description of the cusp groups $\Gamma_\infty(d)$ of the Picard modular groups with $d=1,3$ above, we get:
	
	\begin{cor}\label{NAcusps} The non-arithmetic complex triangle groups with $p=6$ have a cusp group isomorphic to $\Gamma_\infty(3)=G^0_\infty(\mathcal{O}_3)$. The non-arithmetic complex triangle groups with $p=4$ have a cusp group containing an index-2 subgroup isomorphic to $\Gamma_\infty(1)=G^0_\infty(\mathcal{O}_1)$.
	\end{cor}
	
	In order to realize all Nil-manifold groups as cusp groups of non-arithmetic orbifolds in Section~\ref{NACusps} we will also need to use the fact that these non-arithmetic lattices are \emph{generated by $\R$-reflections}, in the sense that they have index 2 in a group generated by $\R$-reflections (recall that the latter are the anitholomorphic involutions in ${\rm Isom}({\rm} H^2_\C)$). More specifically, it was shown in Section~2.2  of \cite{DFP} that any \emph{symmetric complex triangle group} $\Gamma(p,t)$ or $\mathcal{S}(p,\tau)$ is contained (with index 2 or 6) in a group $\widetilde{\Gamma}(p,t)$ or $\widetilde{\mathcal{S}}(p,\tau)$  generated by (three) $\R$-reflections. Moreover one of these three $\R$-reflections, $\sigma_{12}$ in the notation of \cite{DFP}, stabilizes the cusp point $q_\infty$. From this (using again Table~11 of \cite{DPP2}) we record the following:
	
	\begin{cor}\label{NAcusps2} 
	\begin{itemize}
	\item The non-arithmetic lattices $\widetilde{\Gamma}(6,1/6)$, $\widetilde{\mathcal{S}}(6,\sigma_1)$ and $\widetilde{\mathcal{S}}(6,\overline{\sigma_4})$ have a cusp group isomorphic to $G_\infty(\mathcal{O}_3)$. 
	\item The non-arithmetic lattices $\widetilde{\mathcal{S}}(4,\sigma_1)$, $\widetilde{\mathcal{S}}(4,\sigma_5)$ and $\widetilde{\mathcal{S}}(4,\overline{\sigma_4})$ have a cusp group isomorphic to $G_\infty(\mathcal{O}_1)$.
	\end{itemize}
	\end{cor}
	﻿

	\section{Nil-manifolds and cusp cross-sections}
	\label{sec:results}
	
	In this section we prove Theorems \ref{1234} and \ref{567}. We start by recalling some relevant results from \cite{M1} and \cite{M2}. In the context of complex hyperbolic 2-orbifolds/manifolds and their cusp cross-sections, the main result of \cite{M1} and \cite{M2} is the following (Corollary 1.6 of \cite{M2}).
	﻿
	\begin{thm}[\cite{M2}]\label{McRMain} Every Nil-manifold is diffeomorphic to to a cusp cross-section of an arithmetic complex hyperbolic 2-manifold.
	\end{thm}
	﻿
	﻿
	The proof of Theorem \ref{McRMain} proceeds in three main steps: (1) Given any Nil-manifold $N$, find a faithful representation $\rho : \pi_1(N) \longrightarrow \Gamma$ for some arithmetic lattice $\Gamma$ in ${\rm PU}(2,1)$. Geometrically, this realizes a $\pi_1$-injective immersion of $N$ into a cusp cross-section of the orbifold $M={\rm H}^2_\C/\Gamma$. Step (2) is then to show that this immersion can be promoted to an embedding in a finite cover $M'={\rm H}^2_\C/\Gamma'$ of $M$ using a separability argument, so that $N$ is now embedded as the cross-section.of a cusp of $M'$. Step (3), accomplished in \cite{M2}, is to refine the previous step to a further \emph{manifold} cover $M''$, so that $N$ is a cusp cross-section of the manifold $M''$.
	﻿
	
	We now give more details about step (1) using some notation introduced earlier, as we will be extending the results from this step in the next section. 
	Denote as above $G^0_\infty(\mathcal{O}_d) = G^0_\infty \cap {\rm PU}(2,1,\mathcal{O}_d)$ and $G_\infty(\mathcal{O}_d) = \langle G^0_\infty(\mathcal{O}_d) ,\sigma \rangle$, with $G_\infty, G^0_\infty$ defined in (\ref{ParabolicSubgroups}). 
	Note that, from the matrix form~\ref{Pztu}, $G^0_\infty(\mathcal{O}_d)=\{ (z,t,u) \, | \, z,(-|z|^2+it)/2, u \in \mathcal{O}_d \}$. 
	﻿
	
	Recall that a \emph{holonomy representation} of a Nil-manifold group $\pi_1(N)$ is a representation $\rho: \pi_1(N) \longrightarrow G_\infty \simeq {\rm Isom}({\rm Nil})$ such that ${\rm Nil}/\rho(\pi_1(N))$ is diffeomorphic to $N$. By Lemma~\ref{holonomy}, this is equivalent to $\rho$ being discrete, faithful and cocompact.
	﻿
	
	The first part of the following statement is Theorem~7.1 of \cite{M1}. The second part collects case-by-case information from the explicit representations given there in the Appendix.
	\begin{thm}[\cite{M1}]\label{McRFaithful} For any Nil-manifold $N$, there exists a faithful (discrete, cocompact) representation $\rho: \pi_1(N) \longrightarrow G_\infty(\mathcal{O}_d)$ with $d=1$ or 3. More precisely, when $N$ belongs to families (1), (2) or (5) (Nil-torus, vertical half-twist, or $1/4$-twist) the image of $\rho$ is contained in $G_\infty^0(\mathcal{O}_1)$; when $N$ belongs to families (3) or (4) (horizontal half-twist or double half-twist) the image of $\rho$ is contained in $G_\infty(\mathcal{O}_1)$, and when $N$ belongs to families  (6) or (7) ($1/3$-twist or $1/6$-twist) the image of $\rho$ is contained in $G_\infty^0(\mathcal{O}_3)$.
	\end{thm}
	
	\subsection{Nil-manifolds and commensurability classes of Picard modular groups}
	﻿
	We are now ready to prove our main results.
	
	\begin{repthm}{1234}
		\onetwothreefour
	\end{repthm}
	
	In fact we show the following, which implies the conclusion of the theorem. This is a generalization of Theorem~\ref{McRFaithful} for arbitrary $d$, and Steps (2) and (3) described above as in \cite{M1} and \cite{M2} imply that $N$ is then diffeomorphic to a cusp cross-section of a manifold in the given commensurability class.  
	﻿
	﻿
	\begin{lem} The Nil-torus and vertical half-twist groups admit a holonomy representation into $G^0_\infty(\mathcal{O}_d)$ for any $d$; the horizontal and double half-twist groups admit a holonomy representation into $G_\infty(\mathcal{O}_d)$ for any $d$. More precisely, with the notation of the presentations in Section~\ref{nil}, these are given by:
		\begin{enumerate}
			\item (Nil-torus) $a \mapsto (2k,0,1)$, $b \mapsto (2ki\sqrt{d},0,1)$, $c \mapsto (0,16k\sqrt{d},1)$,
			\item (Vertical half-twist) $a \mapsto (2k,0,1)$, $b \mapsto (2ki\sqrt{d},0,1)$, $c \mapsto (0,16k\sqrt{d},1)$, $\alpha \mapsto (0,8k\sqrt{d},-1)$,
			\item (Horizontal half-twist) $a \mapsto (2k,0,1)$, $b \mapsto (2ki\sqrt{d},0,1)$, $c \mapsto (0,8k\sqrt{d},1)$, $\alpha \mapsto (k,0,1)\sigma$,
			\item (Double half-twist) $a \mapsto (2k,0,1)$, $b \mapsto (2ki\sqrt{d},0,1)$, $c \mapsto (0,8k\sqrt{d},1)$, $\alpha \mapsto (-k(1+i\sqrt{d}),4k\sqrt{d},-1)$, $\beta \mapsto (k,0,1)\sigma$.
		\end{enumerate}
	\end{lem}			
	﻿
	\Pf The fact that these formulas define representations, that is, that the images of the generators satisfy the relations of the presentation, can directly be checked using Lemma~\ref{ParabolicProducts}. We give more details below for the double half-turn family (4), in particular since the representation given in \cite{M1} for that family is not correct (the matrix given there in the Appendix for $\alpha$ does not satisfy the relation $\alpha^2=c$).
	﻿
	
	The fact that these representations are faithful follows directly from Lemma 7.2 of \cite{M1}. They are clearly discrete as their image belongs to the discrete subgroup $G_\infty(\mathcal{O}_d)$, and they have cocompact image in Nil as the images of $a,b,c$ already generate a lattice in Nil. Therefore, by Lemma~\ref{holonomy} these representations are all holonomy representations of the corresponding $\pi_1(N)$.
	﻿
	
	We now give more details for family (4). We follow essentially the same strategy as in \cite{M1}, by first finding a representation in $G_\infty$, then conjugating it by a suitable dilation to have entries in $\mathcal{O}_d$. We do this in several steps; all computations are reasonably done by hand, using Lemma~\ref{ParabolicProducts}.
	\begin{itemize}
		\item We first take $\rho(a)=(1,0,1)$ and $\rho(b)=(0,i\sqrt{d},1)$. This gives $[\rho(b),\rho(a)]=(0,4\sqrt{d},1)$, so we take $\rho(c)=(0,2\sqrt{d}/k)$ to satisfy the relation $[b,a]=c^{2k}$.
		\item Next we seek $\rho(\alpha)=(z,t,u)$ satisfying the relations: $\alpha^2=c, \alpha a = a^{-1} \alpha c^k \, (1)$, and $\alpha b = b^{-1} \alpha c^{-k} \, (2)$. (The presentation also lists $[c,\alpha]={\rm Id}$ but this follows from $\alpha^2=c$).
		﻿
		The relation $\alpha^2=c$ implies $u = \pm 1$; we take $u=-1$ (or else $\alpha$ would commute with $c$). Taking $t=\sqrt{d}/k$ gives  $\alpha^2=c$ (independently of the value of $z$).
		﻿
		We then solve for $z$ using relations (1) and (2). The former gives ${\rm Im} \, z=-\sqrt{d}/2$ and the latter ${\rm Re} \, z =-1/2$. Therefore, taking $z=-(1+i\sqrt{d})/2$ gives all the above relations.
		﻿
		(Note: the matrix given for $\rho(\alpha)$ in the Appendix of \cite{M1} corresponds to $(z,t,u)=(2k(1+i),0,-1)$, which squares to the identity, not to $\rho(c)=(0,32k,1)$.)
		﻿
		\item We then seek $\rho(\beta)=(z,t,u)\sigma$ satisfying the relations: $\beta^2=a, \beta c = c^{-1} \beta$, and $\beta b = b^{-1} \beta c^{-k}$. (The presentation also lists $[a,\beta]={\rm Id}$ but this follows from $\beta^2=a$).
		﻿
		We try $\beta=(x,0,1)\sigma$ (for some $x \in \R$) as in \cite{M1}. Taking $x=1/2$ gives the relation $\beta^2=a$. We then check that  the other two relations hold.
		﻿
		\item We then check that the last relation $\alpha \beta = a^{-1} b^{-1} \beta \alpha c^{-k-1}$ is also satisfied.
		﻿
		\item Finally, we conjugate by a dilation of factor $2k$ in order for $\rho$ to have image contained in $G_\infty(\mathcal{O}_d)$, giving the formulas in the statement of the Lemma.
	\end{itemize} \EPf		
	
	Note that in the representations above, any Nil-manifold in family (3) or (4) has a generator ($ \alpha $ or $ \beta $, respectively) which maps to an isometry induced by an antiholomorphic isometry of $ {\rm H}^2_\C $. Recall that a  \emph{complex hyperbolic manifold} is a quotient ${\rm H}_\C^n/\Gamma$ where $\Gamma$ is a torsion-free discrete subgroup of the holomorphic isometry group ${\rm PU}(n,1) \simeq {\rm Isom}^0({\rm H}^n_\C)$, and more generally an \emph{${\rm H}^n_\C$-manifold} is a quotient ${\rm H}_\C^n/\Gamma$ where $\Gamma$ is a torsion-free discrete subgroup of ${\rm Isom}({\rm H}^n_\C)$. We slightly abuse terminology by calling holomorphic (resp. antiholomorphic) isometries of Nil which are induced by holomorphic (resp. antiholomorphic) isometries of ${\rm H}^2_\C$.
	﻿
	\begin{repthm}{34antihol}
		\threefour
	\end{repthm}		
	﻿
	\Pf This follows from the second generalized Bieberbach theorem, Theorem~\ref{GenBieberbach} (2) above. Indeed, automorphisms of Nil are induced by isometries of 
	${\rm H}^2_\C$ (see e.g. Section 2.4 of \cite{M1}), hence conjugation by any element of the affine group ${\rm Nil} \rtimes {\rm Aut}({\rm Nil})$ preserves holomorphicity (resp. antiholomorphicity). Moreover, the holonomy representations produced above (and in \cite{M1}) for these two families contain antiholomorphic isometries.
	﻿
	
	We now give a more pedestrian proof using computations similar to the ones above. Auslander showed (Proposition 2 of \cite{A}) that, if $\Gamma$ is a discrete cocompact subgroup of $N \rtimes C$ (with as previously, $N$ a simply-connected nilpotent Lie group and $C$ a compact Lie subgroup of ${\rm Aut}(N)$), then $\Gamma \cap N$ is a maximal nilpotent normal subgroup of $\Gamma$. (Dekimpe points out in \cite{Dek} that one can remove the word ``normal", with the same proof). In the notation of the presentations of Nil-manifold groups given in Section~\ref{nil}, the maximal nilpotent subgroup (called \emph{fitting subgroup} in \cite{Dek} and \cite{M1}) is $\langle a,b,c \rangle$.
	﻿
	
	By this result, if $\rho: \pi_1(N) \longrightarrow G_\infty$ is a holonomy representation of $\pi_1(N)$ (that is, by Lemma~\ref{holonomy}, discrete, faithful and cocompact) then $\rho(a),\rho(b),\rho(c)$ must be Heisenberg translations. Moreover, $\rho(a)$ and $\rho(b)$ must be horizontal Heisenberg translations, or else they would commute with each other (and $\rho$ would not be faithful). That is, they are of the form $(z,t,1)$ with $z \neq 0$ in $(z,t,u)$-coordinates on ${\rm Nil} \rtimes {\rm U}(1)$. Now, in the horizontal half-twist family (3), the element $\alpha$ satisfies $\alpha^2=a$. Assuming that $\rho(\alpha)$ is a holomorphic isometry of Nil, that is of the form $(z,t,u)\in {\rm Nil} \rtimes {\rm U}(1)$, by Lemma~\ref{ParabolicProducts} we would have $\rho(\alpha)^2=(uz+z,2t,u^2)$. Therefore $u^2=1$, that is $u=\pm 1$. But if $u=1$, $\rho(\alpha)$ would commute with $\rho(c)$, contradicting again faithfulness of $\rho$ (in the presentation, $\alpha$ and $c$ satisfy $\alpha c =c^{-1} \alpha$). Therefore $u=-1$, hence $\rho(\alpha)^2=(0,2t,1)$, a vertical Heisenberg translation, contradicting the fact that $\rho(a)$ is a horizontal Heisenberg translation. Therefore $\rho(\alpha)$ must be antiholomorphic.
	
	﻿
	For the double half-twist family (4) one argues similarly with the element $\rho(\beta)$ satisfying $\rho(\beta)^2=a$. \EPf
	﻿
	
	To complete the proof of Theorem \ref{567}, we must obstruct certain Nil-manifolds from appearing as cross-sections in all but one commensurability class (of Picard modular groups).
	Note that it is not enough to obstruct them from appearing as a subgroup of the ``standard cusp" (represented above by the vector $(1,0,0)^T$); the standard cusp is the only cusp of $X_d ={\rm H}_\C^2 / {\rm PU}(2, 1, \mathcal{O}_d)$ exactly when the ideal class number of $ \mathbb{Q}(i\sqrt{d}) $ is 1 (that is, when $ d = 1, 2, 3, 7, 11, 19, 43, 67, $ or $ 163 $).

	\begin{repthm}{567}
		\fivesixseven
	\end{repthm}
	\Pf The forward direction - that the Nil-manifolds do occur in the single commensurability class given - was proven directly in \cite{M1}, providing explicit holonomy representations into the corresponding $G^0_\infty(\mathcal{O}_d) = G^0_\infty \cap {\rm PU}(2,1,\mathcal{O}_d)$ (Section 9 of \cite{M1}).  We provide explicit representations here, because some of the ones given in \cite{M1} had minor errors, and because they were given in a different form (as matrices in $ U(2,1) $ rather than as points in $ \text{Isom}(\text{Nil}) $).  With the notation of the presentations in Section \ref{nil}, we have:
	
	\begin{enumerate}
		\addtocounter{enumi}{4}
		\item ($1/4$-twist) $a \mapsto (2k,0,1)$, $b \mapsto (2ki,0,1)$, $c \mapsto (0,16k,1)$, $\alpha \mapsto (0,4pk,i)$,
		\item ($1/3$-twist) $a \mapsto (24k,0,1)$, $b \mapsto ((-12 + 12\sqrt{-3})k,0,1)$, $c \mapsto (0,1152k\sqrt{3},1)$, \\ $\alpha \mapsto (-6k + 6k\sqrt{-3} - 12k_1 - 12k_1\sqrt{-3}, 384kk_2\sqrt{3} + 48k^2\sqrt{3} - 96kk_1\sqrt{3} + 192k_1^2\sqrt{3}, \frac{-1+\sqrt{-3}}{2})$,
		\item ($1/6$-twist) $a \mapsto (12k,0,1)$, $b \mapsto ((-6 + 6\sqrt{-3})k,0,1)$, $c \mapsto (0,288k\sqrt{3},1)$, \\ $\alpha \mapsto (6k, 36k^2\sqrt{3} + 48kk_1 \sqrt{3}, \frac{1+\sqrt{-3}}{2}) $.
	\end{enumerate}
	
	
	Now let $N$ be a Nil-manifold in family (5), (6) or (7), and let $\rho: \pi_1(N) \longrightarrow G_\infty \simeq {\rm Isom}({\rm Nil})$ 
	be a holonomy representation of $\pi_1(N)$ with image contained in a subgroup $\Gamma$ of ${\rm Isom}({\rm H}^2_\C)$ commensurable to a Picard modular group ${\rm PU}(2, 1, \mathcal{O}_d)$. As in the proof of the previous theorem, by the second generalized Bieberbach theorem all elements of $\rho(\pi_1(N))$ must be holomorphic, since this is the case for the holonomy representations found in \cite{M1}. Therefore we may assume that $\Gamma$ is contained in PU(2,1); since it is also commensurable to  ${\rm PU}(2, 1, \mathcal{O}_d)$, $\Gamma$ is contained in ${\rm PU}(2, 1, E_d)$ where $E_d=\Q(i\sqrt{d})$ is the field of definition of the commensurability class (see Section 3.1 of \cite{St2}).
	﻿
	
	Now let $A \in {\rm U}(2, 1, E_d)$ be the matrix representative of $\rho(\alpha)$ of the form $P_{(z,t,u)}$ as in Equation~\ref{Pztu}, with $(z,t,u) \in {\rm Nil} \rtimes U(1)$. Again by the second generalized Bieberbach theorem, $\rho(\alpha)$ is affinely conjugate to the corresponding $P_{(z,t,u)}$ for the relevant McReynolds representation, which as noted above has  $u=\zeta_4, \zeta_3,\zeta_6$ for the $1/4$-twist,  $1/3$-twist and $1/6$-twist families respectively. Note that the matrix $P_{(z,t,u)}$ has eigenvalues $1,1,u$. Since, as also noted above, automorphisms of Nil are induced by isometries of ${\rm H}^2_\C$ (which are in turn either matrices or matrices preceded by complex conjugation of coordinates), affine conjugation by an element of ${\rm Nil} \rtimes {\rm Aut} ({\rm Nil})$ either preserves the eigenvalues of the matrix $A$, or maps them to their complex conjugates. Therefore the matrix $A$ has eigenvalues $1,1,u$ with $u=\zeta_4, \zeta_3,\zeta_6$ or their complex conjugates.
	﻿
	
	But the eigenvalues of $A$ must lie in a cubic extension of $E_d$, since $A$ is a $3 \times 3$ matrix with entries in $ E_d $. However the elements $u=\zeta_4, \zeta_3,\zeta_6$ and their complex conjugates do not lie in any cubic extension of $E_d$, except $E_4=\Q(i)$ itself for $\zeta_4$ and its complex conjugate, and $E_3=\Q(\zeta_3)$ for $\zeta_3,\zeta_6$ and their complex conjugates. Therefore $d=1$ if $N$ belongs to the $1/4$-twist family, and $d=3$  if $N$ belongs to the $1/3$- or $1/6$-twist family. \EPf
	﻿
	\subsection{Nil-manifolds and cusp cross-sections of non-arithmetic manifolds}\label{NACusps}
	﻿
	We finish with the following result about cusps of non-arithmetic manifolds.
	﻿
	\begin{repthm}{NA}
		\nonarith
	\end{repthm}
	﻿
	\Pf This follows the strategy of proof of Theorem~\ref{McRMain} given in \cite{M1} and \cite{M2}, outlined after the statement of Theorem~\ref{McRMain}. Step (1) essentially follows from Theorem~\ref{McRFaithful} and Corollary~\ref{NAcusps}. More precisely, by  Theorem~\ref{McRFaithful}, for any Nil-manifold $N$ there exists a faithful (discrete, cocompact) representation $\rho_1: \pi_1(N) \longrightarrow G_\infty(\mathcal{O}_d)$ with $d=1$ or 3. But, by Corollary~\ref{NAcusps2}, there exist non-arithmetic lattices $\Gamma$ in ${\rm PU}(2,1)$ whose cusp group is exactly $G_\infty(\mathcal{O}_d)$ with $d=1$ or 3. This provides a discrete, faithful representation $\rho_2: \pi_1(N) \longrightarrow \Gamma_\infty$ with $\Gamma_\infty$ a cusp subgroup of a non-arithmetic lattice $\Gamma$.

We want to promote this, as in steps (2) and (3) the proof of Theorem~\ref{McRMain}, to show that $N$ occurs as a cross-section of a (non-arithmetic) \emph{manifold} covering the non-arithmetic orbifold ${\rm H}^2_\C/\Gamma$. Step (2) is to find a finite-index subgroup $\Gamma_1 < \Gamma$ (possibly with torsion) such that $\pi_1(N)$ is the full cusp subgroup of $\Gamma_1$, promoting the immersion of $N$ to an embedding as in Theorem~3.12 of \cite{M1}. Step (3) is to find a further finite-index subgroup $\Gamma_2 \leqslant \Gamma_1$ which is torsion-free but still has a cusp subgroup isomorphic to $\pi_1(N)$, as in Corollary~1.6 of \cite{M2}. 
	
	We now give more details for these last two steps, as the relevant results of \cite{M1} and \cite{M2} are stated there for arithmetic lattices. For step (2) we use the following result of Bergeron (``Lemme principal" of \cite{B}, as stated in Proposition 3.8 of \cite{M1}):
	
	\begin{prop}[\cite{B}] Let $H$ be an algebraic subgroup of a linear algebraic group $G$ and $\Gamma$ a finitely generated subgroup of $G$. Then $H \cap \Gamma$ is separable in $\Gamma$.
	\end{prop}
	Taking $G={\rm SU(2,1)}$, $H = {\rm Stab}_{G}(q_\infty)$ and $\Gamma$ any lattice in $G$ (arithmetic or not), this gives that the cusp subgroup $\Gamma_\infty=\Gamma \cap H$ is separable in $\Gamma$. (Note that $H$ is an algebraic subgroup of $G$ as it consists of all upper-triangular matrices in $G$). The result of step (2) then follows; that is, Theorem~3.12 of \cite{M1} holds without the arithmetic assumption, with the same proof. We now have a non-arithmetic lattice $\Gamma_1$ in SU(2,1) such that $N$ is diffeomorphic to a cusp cross-section of the orbifold ${\rm H}^2_\C/\Gamma_1$. In other words, we have a discrete, faithful representation $\rho_3: \pi_1(N) \longrightarrow \Gamma_1$ such that $\rho_3(\pi_1(N))$ is exactly a cusp stabilizer in $\Gamma_1$.

	For Step (3) we show that the relevant results of \cite{M2} still apply in the non-arithmetic but \emph{integral} case. We first recall the following result, Corollary~5.3 of \cite{M2}:
	
	\begin{thm}[\cite{M2}]\label{manifoldcover} Let $\Lambda_1$ be a torsion-free virtually unipotent subgroup of an arithmetic lattice $\Lambda$. Then there exists a torsion-free finite-index subgroup $\Lambda_0$ of $\Lambda$ containing $\Lambda_1$.
	\end{thm}
	
	We now recall that all known non-arithmetic lattices in SU(2,1) are \emph{integral} in the following sense (see e.g. Section~3.1 of \cite{DPP1}). For each such lattice $\Gamma_1$ there exists a totally real number field $k$ (in fact, the adjoint trace field $\Q[{\rm Tr} \, {\rm Ad} \, \Gamma_1]$) and an imaginary quadratic extension $L$ of $k$ such that $\Gamma_1 < {\rm SU}(H,\mathcal{O}_L)$ for some Hermitian form $H$ of signature (2,1) defined over $L$. Denoting ${\rm Gal}(L/\Q) = \{{\rm Id}=\sigma_1,...,\sigma_n \}$, this provides an injective homomorphism:
$$	
\begin{array}{rrcl} \varphi : & \Gamma_1 & \longrightarrow & \Pi_{i=1}^n {\rm SU}(^{\sigma_i} H, \mathcal{O}_L) \\
 & \gamma & \longmapsto & (^{\sigma_1}\gamma,...,^{\sigma_n}\gamma )
\end{array}
$$	
where, for any $\sigma \in {\rm Gal}(L/\Q)$ and $\gamma \in \Gamma_1$, $^{\sigma} \gamma$ denotes the matrix obtained from $\gamma$ by applying $\sigma$ to the entries of $\gamma$; likewise $^{\sigma} H$ denotes the Hermitian matrix obtained from $H$ by applying $\sigma$ to the entries of $H$ (this is still a Hermitian matrix by the assumptions on $L$, as each $\sigma$ commutes with complex conjugation).

Now $\Lambda = \Pi_{i=1}^n {\rm SU}(^{\sigma_i} H, \mathcal{O}_L)$ is an arithmetic lattice. Indeed, by restriction of scalars (more specifically, choosing a $\Q$-basis for $L$ which is also a $\Z$-basis for $\mathcal{O}_L$), we get a linear algebraic group $G'$ defined over $\Q$ such that $\Lambda \simeq G'(\Z)$, which is therefore an arithmetic lattice in $G'(\R) \simeq  \Pi_{i=1}^n {\rm SU}(^{\sigma_i} H)$ by the Borel--Harish-Chandra theorem (see e.g. Section~5.5 of \cite{WM} for this now classical construction). Note that $\varphi (\Gamma_1)$ has infinite index in $\Lambda$ since $\Gamma_1$ is non-arithmetic.

Note also that for any $\sigma \in {\rm Gal}(L/\Q)$ and $\gamma \in \Gamma_1$, $^{\sigma} \gamma$ is unipotent (resp. finite-order) if and only $\gamma$ is, hence $\varphi(\gamma)$ is unipotent (resp. finite-order) if and only $\gamma$ is. Therefore $\Lambda_1 = \varphi(\rho_3(\pi_1(N)))$ is a torsion-free, virtually unipotent subgroup of the arithmetic lattice $\Lambda$. By Theorem~\ref{manifoldcover}, there exists a torsion-free finite-index subgroup $\Lambda_0$ of $\Lambda$ containing $\Lambda_1$. Taking the intersection with the first factor of $\Lambda$ gives a torsion-free finite-index subgroup $\Gamma_2$ of $\Gamma_1$, containing $\rho_3(\pi_1(N))$ as required for the conclusion of Step (3). \EPf
	
	\raggedright
	
	\vspace{.2cm}
	﻿
	\begin{flushleft}
		\textsc{Julien Paupert\\
			School of Mathematical and Statistical Sciences, Arizona State University}\\
		\verb|paupert@asu.edu|
	\end{flushleft}
	﻿
	\vspace{.2cm}
	﻿
	\begin{flushleft}
		\textsc{Connor Sell\\
			D\'epartement de Math\'ematiques, Universit\'e du Qu\'ebec \`a Montr\'eal}\\
		\verb|sell.connor_daniel@courrier.uqam.ca| 
	\end{flushleft}
	﻿
	
\end{document}